\documentclass[11pt]{article}
\usepackage[margin=4.46cm]{geometry}
\usepackage{graphicx}
\usepackage{epsfig}
\usepackage{amsmath}
\usepackage{amssymb}
\newcommand{\commentout}[1]{}

\def \Rset {{\mathbb R}}

\def \Zset {{\mathbb Z}}

\def \Nset {{\mathbb N}}

\newcommand{\ra}{\rightarrow}

\newcommand{\nit}{\noindent}
\newcommand{\no}{\nonumber}
\newcommand{\be}{\begin{equation}}
\newcommand{\ee}{\end{equation}}
\newcommand{\ba}{\begin{eqnarray}}
\newcommand{\ea}{\end{eqnarray}}
\newcommand{\bi}{\begin{itemize}}
\newcommand{\ei}{\end{itemize}}
\newcommand{\br}{\begin{eqnarray}}
\newcommand{\er}{\end{eqnarray}}

\newcommand{\eps}{\mbox{$\epsilon$}}

\newcommand{\qed}{\mbox{$\square$}\newline}

\newcommand{\real}{\mathbb{R}}

\newtheorem{example}{Example}[section]
\newtheorem{theo}{Theorem}[section]
\newtheorem{defin}{Definition}[section]

\newtheorem{lem}{Lemma}[section]
\newtheorem{cor}{Corollary}[section]
\newtheorem{rmk}{Remark}[section]

\begin{document}
\title{Analysis and Comparison of Large Time Front Speeds in Turbulent Combustion Models}
\author{Jack Xin$^{*}$ \; and \; Yifeng Yu \thanks{Department of Mathematics,
UC Irvine, Irvine, CA 92697, USA. Email: jxin@math.uci.edu, yyu1@math.uci.edu.}}
\date{}
\maketitle
\begin{abstract}
Predicting turbulent flame speed (the large time front speed) is a fundamental problem in
turbulent combustion theory. Several models have been proposed to study the turbulent
flame speed, such as the G-equations, the F-equations (Majda-Souganidis model) and
reaction-diffusion-advection (RDA) equations. G-equations describe the front motion law in
the form of local normal velocity equal to a constant (laminar speed) plus the
normal projection of fluid velocity. In level set formulation, G-equations are
Hamilton-Jacobi equations with convex ($L^1$ type) but non-coercive Hamiltonians.
In the first part of this paper, we show that flow induced strain
reduces front speeds of G-equations in periodic compressible and shear flows.
The F-equations arise in asymptotic analysis of reaction-diffusion-advection equations
and are quadratically nonlinear analogues of the G-equations.  In the second part of the paper,
we compare asymptotic growth rates of the turbulent flame speeds from the G-equations,
the F-equations and the RDA equations in the large amplitude ($A$) regime of
spatially periodic flows. The F and G equations share the same asymptotic front speed growth rate; in particular,
the same sublinear growth law $A\over \log(A)$ holds in cellular flows. Moreover, in two space dimensions,
if one of these three models (G-equation, F-equation and the RDA equation) predicts the bending
effect (sublinear growth in the large flow), so will the other two. The nonoccurrence
of speed bending is characterized by the existence of periodic orbits on the torus and
the property of their rotation vectors in the advective flow fields.
The cat's eye flow is discussed as
a typical example of directional dependence of the front speed bending.
The large time front speeds of the viscous F-equation
have the same growth rate as those of the inviscid F and G-equations in two dimensional
periodic incompressible flows.
\end{abstract}
\thispagestyle{empty}
\newpage
\section{Introduction}
\setcounter{equation}{0}
\setcounter{page}{1}
Turbulent combustion is a complex nonlinear and multiscale phenomenon \cite{Pet00}. A comprehensive
physical-chemical modeling requires a system of reaction-diffusion-advection (RDA) equations coupled with the Navier-Stokes equations.
For theoretical understanding and efficient modeling of the turbulent flame propagation,
various simplified or phenomenological models have been proposed and studied. Most notably, these models
are passive scalar reaction-diffusion-advection equations (RDA) and Hamilton-Jacobi equations (HJ), as documented in
books \cite{W85,Pet00,Xin_09} and research papers \cite{Ab_02,ABP_00,CW,Const,EMS,MS2,NX_09a,OF,R,Siv,S89,Yak} to
name a few. For simplicity of presentation,
throughout this paper, we shall consider the fluid velocity $V$ as
time-independent $V=V(x)$.
\medskip

{\nit} $\bullet$ ({\bf Scalar RDA model}).  The passive scalar reaction-diffusion advection equation for the temperature field is:
\ba
& & T_t +  V(x)\cdot D T= d\, \Delta T+\, {1\over \tau_r}\, f(T), \no \\
& & T(x,0) = T_0(x), \;\; x\in \real^n, \label{rd1}
\ea
where $T$ represents the reactant temperature, $D$ is the spatial gradient operator, $V(x)$ is a prescribed fluid velocity, $f$ is a
nonlinear reaction  function; $d$ is the molecular diffusion constant, $\tau_r >0$ is
reaction time scale. The flow field $V$ is known or statistically known.
For an isothermal reaction, the scalar is a reactant concentration
however we shall still use $T$ to denote it. The common form of the reaction function
is $f(T)=T(1-T)$, so called Kolmogorov-Petrovsky-Piskunov-Fisher (KPP-Fisher);
$f(T)=T^m(1-T)$ ($m\geq 2$, higher order KPP-Fisher); $f(T)=e^{-E/T}(1-T)$ ($E>0$), Arrhenius combustion
nonlinearity; $F(T)=0$, $T\in [0,\theta]\cup \{1\}$, $f(T) > 0$, $T\in (\theta, 1)$, ignition
combustion nonlinearity. KPP or generalized KPP comes from isothermal autocatalytic reaction-diffusion system with
equal diffusion constants (or unit Lewis number), \cite{Bill_91,Xin_09}. Equation (\ref{rd1}) is well-known to admit
propagating front solutions if the advection is absent ($V=0$), \cite{AW_75}.
This is the case when a car engine filled with liquid or gaseous fuel
is ignited to generate a spreading flame. Turbulent combustion concerns with the setting of flame propagation
when the fuel is stirred on a broad range of scales for the purpose of speed enhancement and
waste gas reduction. Though the flame front will be wrinkled by the fluid velocity, its average location
eventually moves at a steady speed $s_T$ in each specified direction, the so called ``turbulent flame speed".
The prediction of the turbulent flame speed is a fundamental problem in turbulent combustion theory \cite{W85,R,Pet00}.
For KPP nonlinearity, it is known \cite{GrFr,BH07,NX_09a,Xin_09} that $s_T$ is given by a variational principle
on the large time growth rate of a viscous quadratically nonlinear Hamilton-Jacobi equation (QHJ). More precisely,
considering compactly supported initial data $T(x,0)$, then for
each direction $e$ and wave number $\lambda > 0$, let $\mathcal{L}_e(\lambda)$ be the principal Lyapunov exponent of
the linear advection-diffusion equation:
$$
\phi_t = d\, \Delta \phi +(-2\, d\, \lambda \, e - V(x,t))\cdot D \phi+ [d\, \lambda^2 + \lambda \,
e\cdot  V(x,y) + \tau^{-1}_{r} f' (0)]\, \phi,
$$
with initial data $\phi(x,0)=1$. The function $u = \ln \phi$ then satisfies the viscous QHJ:
$$
u_t = d\, \Delta u + d\, |D u |^2 + (-2\, d\, \lambda \, e - V(x,t))\cdot Du +
d\, \lambda^2 + \lambda \, e\cdot  V(x,y) + \tau_{r}^{-1}\, f'(0).
$$
Under suitable stationarity and ergodicity condition of the
flow field \cite{NX_09a}, the following limit exists almost surely and is independent of $x$:
$$
\mathcal{L}_e(\lambda) = \lim_{t \ra +\infty} {1\over t}u(x,t).
$$
The turbulent front speed along in $e$ is:
\be
s_T (e) = \inf_{\lambda > 0} {\mathcal{L}_e(\lambda) \over \lambda}, \label{rd2}
\ee
a deterministic quantity at large scale. In particular, if $V$ is a periodic flow field, then
$$
\mathcal {L}_{e}(\lambda)=\tau_{r}^{-1}\, f'(0)+H^{*} (\lambda e).
$$
Here for $p\in \Rset^n$,  $H^{*}(p)$ is given by the following cell problem
$$
-d\Delta w+d|p+Dw|^2+V(x)\cdot (p+Dw)=H^{*}(p)
$$
for $w\in C^{\infty}(\Bbb T^n)$, where $\Bbb T^n$ is the $n$-dimensional flat torus.
\medskip

\nit $\bullet$ ({\bf Majda-Sougandis model}).  Turbulent combustion might involve
small turbulent scales which are no less than the reaction time scale in the thin flame model.
In \cite{MS1}, the velocity field $V$ is space-time periodic and
in a scale-separation form
$V=V(x,t,\eps^{-\alpha} x, \eps^{-\alpha} t)$, $d= \eps\, d$, $\tau_r = \eps \tau_r$,
$\alpha \in (0,1]$.  For simplicity of analyzing upscaling,  we consider
$V$ as time-independent with no integral scales, i.e $V=V(\eps^{-\alpha} x)$.
The limiting behavior of $T=T^{\eps}$ is \cite{MS1}:
$\lim_{\eps \ra 0} T^{\eps} = 0$ locally uniformly in $\{(x,t):Z < 0\}$
and $T^{\eps} \ra 1$ locally uniformly in the interior of
$\{(x,t): Z=0\}$, where $Z \in C(\real^n\times [0,+\infty))$ is the unique
viscosity solution of the variational inequality
\be
 \max (Z_{t} - \bar{H}(D Z) -\tau_{r}^{-1}f'(0), Z) = 0, \quad (x,t) \times \real^n\times
 (0,+\infty), \label{ve}
\ee
with initial data $Z(x,0) = 0$ in the support of $T(x,0)$, and $Z(x,0)=-\infty$ otherwise.
The set $\Gamma_{t} =\partial \{ x \in \real^n: Z(x,t) < 0\}$ can be regarded as
a front which moves with normal velocity
$$
v_{n}={\mathcal F}(n)
$$
where $n$ denoting the normal vector pointing to the propagation direction and
$$
{\mathcal F}(p)=\inf_{\lambda>0}{\tau_{r}^{-1}f'(0)+\hat H (p\lambda)\over \lambda}.
$$
The effective Hamiltonian
$\hat {H} = \hat {H}(p)$ is defined as a solution of
the following cell problem: for each $p\in \Rset^n$,  there are a unique number $\hat H(p)$ and a function
$F(x) \in C^{0,1}(\Bbb T^n)$ such that
\be
- a(\alpha)d\Delta F+d|p + DF|^{2} +
V(x)\cdot (p + DF) = \hat H(p), \label{fe}
\ee
where $a(\alpha ) = 0$ if $\alpha \in (0,1)$, $a(\alpha = 1) =1$. The turbulent flame speed along the unit direction $e$ is:
\be\label{MSturbulent}
s_T(e)=\mathcal {F}(e).
\ee
Note that when $\alpha=1$,  (\ref{MSturbulent}) coincides with (\ref{rd2}).  Hereafter, we shall refer to (\ref{fe}) as
F-equation if $d=0$, and viscous F-equation if $d>0$.
\begin{figure}
\centering
\includegraphics{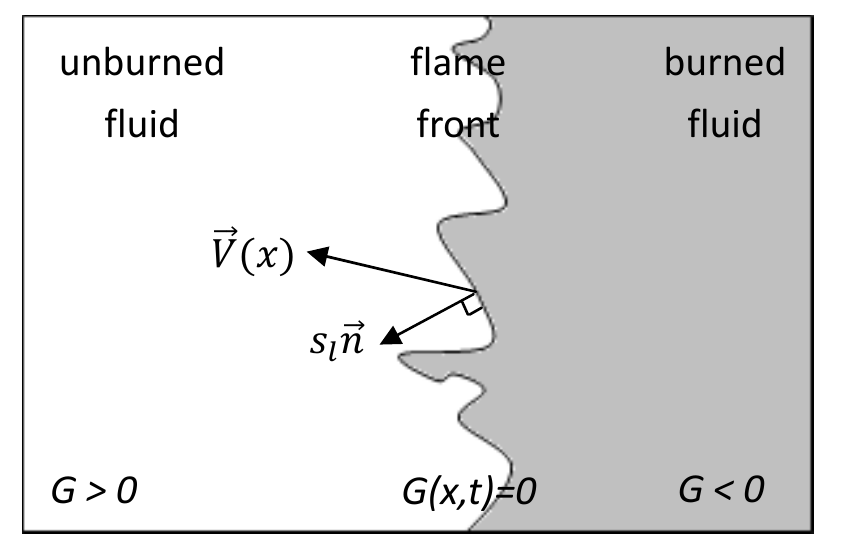}
\caption{Illustration of G-equation (level set) model.}
\label{fig:projection}
\end{figure}
\medskip

\nit {$\bullet$ {\bf G-equation (Level Set) model}}. The RDAs (\ref{rd1}) and in particular the KPP are
first principle equations, yet they are
limited to the unit Lewis number regime.
Another approach in turbulent flame modeling is the level set formulation of interface motion laws in
the thin interface regime. The simplest motion law is that the normal velocity of the interface
($V_n$) is equal to a constant $s_l$ (the laminar speed)
plus the projection of fluid velocity along the normal $\overrightarrow{n}$. See Figure 1. The laminar speed is the flame speed when fluid is
at rest. Let the flame front be the zero level set of a function
$G(x,t)$, the burnt region is $G(x,t) < 0$, and the unburnt region is $G(x,t) > 0$.
The normal direction pointing from the burnt region to the
unburnt region is $DG/|DG|$, the normal velocity is $-G_t/|DG|$.
The motion law becomes the so called $G$-equation, a popular model in turbulent combustion
\cite{W85,Pet00}:
\be
G_t + V(x)\cdot DG + s_l |DG|=0. \label{ge1}
\ee
Chemical kinetics and Lewis number effects are all included in the laminar speed $s_l$ which is
provided by a user. Formally, under the G-equation model, for a specified unit direction $p$,
\be \label{large-T-lim}
s_T(p)=-\lim_{T\to +\infty}{G(x,T)\over T}.
\ee
Here $G(x,t)$ is the solution of equation (\ref{ge1}) with initial data $G(x,0)=p\cdot x$.
We remark that in combustion literature, the level set function is defined as
$G(x,t)> 0$ in the burnt region and $G(x,t)< 0$ in the unburnt region (pp 91-92 of \cite{Pet00}).
The resulting G-equation is $G_t + V(x,t)\cdot DG = s_l |DG|$. We shall work with the form in
(\ref{ge1}) instead.  The existence of $s_T$ has been rigourously established in \cite{XY_10} and \cite{CNS}
independently for incompressible periodic flows. And $s_T$ is the effective Hamiltonian of the following cell problem
\be\label{generalGcell}
s_l|p+DG|+V(x)\cdot (p+DG)=\bar H(p)=s_T.
\ee
Here $\bar H(p)$, the effective Hamiltonian, is the unique number such that the above equation admits periodic approximate solutions. See the Appendix for some basic facts about viscosity solutions and cell problems.
The formal analysis of (\ref{ge1}) and $s_T$ is also performed in the framework of renormalization
group methods \cite{Siv,S89, Yak}. See also \cite{Pet92} on a spectral closure approximation, and
\cite{Bell_08} for a numerical study of G-equation in comparison with a
combustion system modeling thermal-diffusive instabilities of free-propagating
premixed lean hydrogen-air flames. Though G-equation is a phenomenological model, it is more flexible in that
many factors influencing front motion can be incorporated into $s_l$. For example, the strain effect
of a turbulent fluid flow is modeled by extending $s_l$ to $s_l + \, d\overrightarrow{n}\cdot DV\cdot
\overrightarrow{n}$, where $d$ is the Markstein length which is proportional to the flame thickness.
The G-equation with flow induced strain is \cite{Pet00}:
\be\label{strainG}
G_t+s_l|DG|+V(x,t)\cdot DG+\, d{DG\over |DG|}\cdot DV\cdot DG=0.
\ee
Then formally, $s_T$ is also given by (\ref{large-T-lim}), where
$G(x,t)$ is the solution of equation (\ref{strainG}) with initial data $G(x,0)=p\cdot x$.
So far we are not able to prove the existence of $s_T$ except for some simple situations
like the one space dimensional (1d) compressible flow and the incompressible shear flow.
It is conjectured by some experts \cite{R} in combustion theory that
the strain term will slow down flame propagation. Theoretically it is hard to
verify this for general flows. In section 2 on this paper, we confirm this conjecture for 1d compressible flow
and the shear flow. Besides the strain effect,  the flame stretching mechanism
also includes the curvature effect.  One such model proposed in \cite{Pet00} is
to replace $s_l$ by $s_l(1-d\, \kappa)$. Then the G-equation becomes
\be\label{curvatureG}
G_t-ds_l\kappa|DG|+s_l|DG|+V(x,t)\cdot DG=0.
\ee
Here $d$ is the Markstein length and $\kappa$ is the mean curvature of the flame front,
i.e, $\kappa={\rm div}\,({DG\over |DG|})$. The curvature G-equation (\ref{curvatureG}) is very difficult to analyze.  To obtain some ideas of the diffusion effect, a natural simplification is to change the mean curvature term $\kappa$ to $\Delta G$. This leads to the viscous G-equation
\be\label{viscousG}
G_t-ds_l\Delta G+s_l|DG|+V(x,t)\cdot DG=0.
\ee
The above viscous G-equation also serves as a basic model to understand the numerical
diffusion effect introduced in the numerical computation of equation (\ref{ge1}).  For the viscous case,
$s_T=s_T(p,d)$ is given by the large time limit (\ref{large-T-lim}), where
$G(x,t)$ is the solution of equation (\ref{viscousG}) with initial data $G(x,0)=p\cdot x$.
It is also the effective Hamiltonian of the following cell problem
$$
-ds_l\Delta G+s_l|p+DG|+V(x)\cdot (p+DG)=s_T(p,d)=\bar H(p,d).
$$
Here $\bar H(p,d)$ is the unique number such that the above equation admits periodic solutions.
The most general G-equation is to combine both the strain effect and the curvature
effect, that is $s_l$ in the basic G-equation (\ref{ge1}) is replaced by $s_l(1-d\kappa)+\, d\overrightarrow{n}\cdot DV\cdot \overrightarrow{n}$.
\medskip

\medskip

Though various passive scalar models as shown above have been proposed to study $s_T$, their predictions may be
potentially different or sometimes asymptotically identical. It requires delicate analysis to understand these subtleties.
The goal of this paper is to analyze and
compare qualitative and quantitative properties of $s_T$ in G-equations (\ref{ge1})-(\ref{strainG}),
F-equation (\ref{fe})-(\ref{MSturbulent}) and the RDA equation (\ref{rd1})-(\ref{rd2}) for
steady (time-independent) periodic flows ($V=V(x)$) with mean equal to zero.
The $s_T$'s are compared in terms of different nonlinearities, and the flow induced strains.
\medskip

\nit{\bf{Outline of the paper}}. To simplify notations, throughout this paper, we normalize
$s_l=d=\tau_r=1$.
\medskip

\nit{\bf  Section 2:} We analyze $G$-equation (\ref{ge1})-(\ref{strainG}) with and without the strain effect. We
show that the strain term slows down the propagation speed $s_T$
in one-dimensional compressible flows and in shear flows.
\medskip

\nit{\bf Section 3:} The comparison of the F-equation (\ref{fe})-(\ref{MSturbulent}) with the G-equation (\ref{ge1})
had been studied in \cite{EMS} for periodic shear flows. Their computations showed that G-equation model always predicts slower turbulent flame speeds than the F-equation model for the shear flow.   We will prove this result for general flows as long as $f'(0)$ is no less than a threshold value.  We also investigate
the asymptotic behavior of $s_T$ from the reaction-diffusion equation (\ref{rd1})-(\ref{rd2}),
the G-equation (\ref{ge1}) and the F-equation (\ref{fe})-(\ref{MSturbulent})
when $V(x)$ is scaled to $A\, V(x)$, $A \gg 1$.
We prove that the asymptotic growth rate $\lim_{A\to +\infty}{s_T\over A}$ is the same for G and F equations.
The limit is given by an inf-max formula which is also equivalent to a variational formula involving
the invariant measures of the flow $\dot{x}=V(x)$ based on the weak KAM theory.
We also provide a necessary and sufficient condition
for the bending effect (sublinear growth in $A$) in two space dimensions ($n=2$) in terms of the
periodic orbits of the dynamical system $\dot{x}=V(x)$ on the two dimensional
torus and their rotation vectors.  We show that when $n=2$ the asymptotic growth rate of $s_T$ from
the RDA equation is qualitatively the same as $s_T$ from the G and F equations. In particular, these three models will predict the bending effect simultaneously.  See Theorem \ref{character} for the precise statement. As an application, the growth laws of $s_T$  as $a\to +\infty$ predicted by these three models are completely identified for the cat's eye flow and the cellular flow.  In particular, for cellular flows, the ${A\over \log(A)}$ growth law of $s_T$ in G-equation (\ref{ge1}) is also verified for the
F-equation (\ref{fe})-(\ref{MSturbulent}).  These two flows
appeared as representative flow examples in related dynamo and convection-enhanced
diffusion problems \cite{CS_89,fannjiang}.

\medskip
\nit{\bf Section 4:} We prove that in
two space dimensions, the effective Hamiltonian from the viscous F-equation gives the same asymptotic growth
rate as $s_T$ in $A$ predicted by the F-equation (\ref{fe})-(\ref{MSturbulent}) and the (inviscid) G-equation (\ref{ge1}).
It follows that the effective Hamiltonian of the viscous F-equation grows sublinearly in large cellular flows.
On the other hand, it  has a lower bound of $O(A^{1/4})$, see \cite{NR}. Moreover, we also show that, in 2d,
changing the order of  sending the flame thickness or reaction time scale to zero
and sending the turbulent intensity to infinity does not change the asymptotic growth rate of the
front speed $s_T$.

\medskip
\nit{\bf Appendix:} For reader's convenience, we will review some basic facts about
viscosity solutions of Hamilton-Jacobi equations and their cell problems.

\medskip
\nit {\bf Acknowledgements:}
The work was partially supported by NSF grants
DMS-0911277 (JX) and DMS-0901460 (YY).

\medskip

\nit{\bf Assumptions and notations}: Throughout this paper, we assume that nonlinear function $T\to f(T)$ is of KPP type, i.e\\
$$
f\in C^{1,\epsilon}([0,1]),\ f(0)=f(1)=0
$$
and $f$ is non-increasing on $(1-\epsilon,1)$ for some $\epsilon>0$. Also
$$
0<f(s)\leq sf'(0)  \quad \text{for $s\in (0,1)$}.
$$
Note a prototypical example is
$$
f(T)=KT(1-T).
$$
(1) $\Bbb Z^n=\{(x_1,x_2,...,x_n)|\ x_i\in \Zset\}$. $\Bbb T^n$ denotes the $n$-dimensional flat torus, i.e, $\Bbb T^n=\Rset ^n/\Zset^n$.\\
(2) $f\in C^{r}(\Bbb T^n)$ if $f\in C^{r}(\Rset ^n)$ and it is periodic, i.e, $f(x+v)=f(x)$ for $x\in\Rset^n$ and $v\in \Zset^n$.\\
(3) $Df$ is the gradient of $f$. If $g\in \Rset^n\to \Rset$.\\
(4) Throughout this paper, we assume the velocity field $V\in C^{\infty}(\Bbb T^n)$ and has mean zero, i.e, $\int_{\Bbb T^n}V\,dx=0$.\\
(5) For a curve $\xi:\Rset\to \Rset^n$, $\hat \xi$ is its image under the natural
projection $\Rset ^n\to \Bbb T^n$.  See Figure 2 for a projected periodic orbit on the torus.\\
\medskip
\nit{If we scale $V$ to $AV$}:\\
(6) $\alpha_A(p)$ the turbulent flame speed predicted by the G-equation equation model (\ref{ge1})-(\ref{generalGcell}).\\
(7) $\gamma_A(p)$ the turbulent flame speed predicted by the F-equation equation model (\ref{fe})-(\ref{MSturbulent}) for $\alpha\in (0,1)$.\\
(8) $c_{p}^{*}(A)$ the turbulent flame speed predicted by the RDA equation model (\ref{rd1})-(\ref{rd2}).

\begin{figure}
\centering
\includegraphics{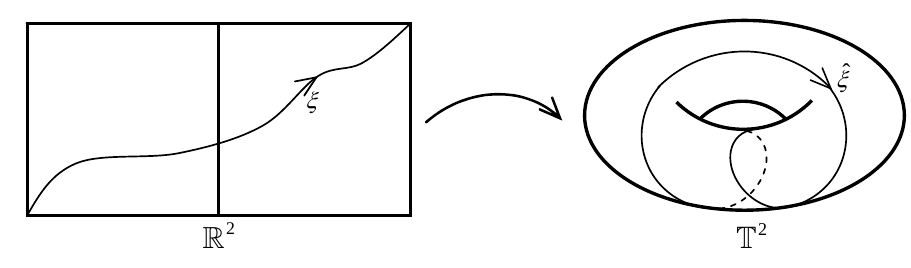}
\caption{Schematic of the projection $\Rset^n\to \Bbb T^n$, and a periodic orbit on $\Bbb T^n$.}
\label{fig:projection}
\end{figure}

\section{G-equations with and without Strain}
\setcounter{equation}{0}
The goal of this section is to show that the strain term in (\ref{strainG}) slows down flame propagation for 1d compressible flow and the shear flow. This is different from the diffusion  effect which actually enhances the flame propagation for 1d compressible flow \cite{LXY}. Recall that the turbulent flame speed predicted by the G-equation is the effective Hamiltonian associated with the cell problem (\ref{generalGcell}). See the Appendix for some basic facts about viscosity solutions and cell problems.
\subsection{One-dimensional G-equation and Compressible Flow}
It suffices to look at $p=1$.  We write the flow velocity $V(x)=v(x)$. For one-dimensional G-equation without the strain term, the turbulent flame speed is the unique number $c$ such that the following equation ($p=1$)
$$
|1+G'|+v(x)(1+G')=c
$$
admits approximate periodic viscosity solutions.
Let us assume that  $\int_{0}^{1}v(x)\,dx=0$ to make $c\geq 0$. It has been proved in section 3 of \cite{LXY} that
$$
c=
\begin{cases}
0  \quad \text{if $\{x|\  v(x)=-1\}\ne \emptyset$}\\
{(\int_{0}^{1}{1\over 1+v(x)}\,dx)}^{-1}>0.
\end{cases}
$$
Now let us consider the G-equation with the strain term.  Using the same method as in section 3 of \cite{LXY},
we can show that there exists a unique number $\hat c$ such that the following equation
$$
(1+v')|1+G'|+v(x)(1+G')=\hat c
$$
admits approximate periodic viscosity solutions. The $\hat c$ is the turbulent flame speed predicted by the G-equation under strain, and is given by
$$
\hat c=
\begin{cases}
0  \quad \text{if $\{x|\  v'(x)+v(x)=-1\}\ne \emptyset$}\\
{(\int_{0}^{1}{1\over 1+v'(x)+v(x)}\,dx)}^{-1}>0.
\end{cases}
$$
\begin{theo}\label{strain1d}
The flame speed in the G-equation with strain is no faster than that of the G-equation without strain:
$$
c\geq \hat c.
$$
If $c>0$, then ``=" holds if and only if $v'\equiv 0$.
\end{theo}
We first prove the following lemma.
\begin{lem}\label{inverse} Suppose that $w\in C^{1}(\Bbb T^1)$ and $w+w'>0$. Then $w>0$ and
$$
\int_{0}^{1}{1\over w}\,dx\leq \int_{0}^{1}{1\over w+w'}\,dx.
$$
The ``=" holds if and only if $w$ is a constant.
\end{lem}
\nit Proof:  Assume that $w(x_0)=\min_{\Bbb T^1}w$.  Then
$$
w(x_0)=w(x_0)+w'(x_0)>0.
$$
So $w>0$.  Now let $s={1\over w}$.  Then ${s^2\over s-s'}={1\over w+w'}$. Since $s-s'>0$,  Cauchy's inequality implies that
$$
\int_{0}^{1}(s-s')\,dx\int_{0}^{1}{s^2\over s-s'}\,dx\geq \left (\int_{0}^{1}s\,dx \right )^2.
$$
Note that $\int_{0}^{1}(s-s')\,dx=\int_{0}^{1}s\,dx$.  Accordingly,
$$
\int_{0}^{1}{s^2\over s-s'}\,dx\geq \int_{0}^{1}s\,dx.
$$
The equality holds if and only there exists $\lambda>0$ such that
$$
s-s'=\lambda s.
$$
Since $s$ is periodic, this implies that $s'\equiv 0$.  So $w'\equiv 0$.
\qed

{\bf \nit Proof of Theorem  \ref{strain1d}} Suppose that $\hat  c>0$. Then $1+v+v'>0$ and the above theorem follows immediately from Lemma \ref{inverse} by setting $w=1+v$. \qed

\subsection{Strain Effects in Shear Flows}
Suppose that $V(x,y)=(v(y),0)$ for $(x,y)\in \Rset^2$. For $(m,n)\in \Rset ^2$,  denote $\lambda (m,n)$ as the unique number such that the G-equation for the shear flow
$$
\sqrt {m^2+(n+u')^2}+mv(x)=\lambda=\lambda (m,n)
$$
has a periodic viscosity solution. We also write $\hat \lambda (m,n)$ as the unique number such that the G-equation for the shear flow with strain term
$$
\sqrt {m^2+(n+u')^2}+{m(n+u')v'\over \sqrt {m^2+(n+u')^2}}+mv(x) =\hat \lambda=\hat \lambda (m,n).
$$
It is clear that
\be\label{lowerbound2}
\lambda,\ \hat \lambda\geq |m|+\max_{\Bbb T^1}mv
\ee
To ensure that both of them are nonnegative, we further assume that
$$
\int_{0}^{1}v(x)\,dx=0.
$$
Note that $\lambda (m,n)$ and $\hat \lambda(m,n)$ are the turbulent flame speeds predicted by the G-equation without and with strain effects respectively for the shear flow.
\begin{theo}\label{shear}
The flame speed under the strain of shear flows is no faster than that in shear flows without the strain.
$$
\lambda (m,n)\geq \hat \lambda (m,n).
$$
If $m\ne0$ and $\lambda (m,n)>|m|+\max_{\Bbb T^1}mv$, then ``=" holds if and only if $v'\equiv 0$.
\end{theo}

We first prove several lemmas.

\begin{lem}\label{unique}
Given $m, c\in \Rset$ and $a>|m|$, there exists a unique $p>0$ such that
$$
\sqrt {m^2+p^2}+{cp\over \sqrt {m^2+p^2}}=a.
$$
\end{lem}
Proof:  Denote $f(p)=\sqrt {m^2+p^2}+{cp\over \sqrt{m^2+p^2}}$.  Clearly there exists at least one $p_0>0$ such that
$$
f(p_0)=a.
$$
Note that if $f(p)=a$ and $p>0$, then $p+c>0$ and
$$
f'(p)={(p+c)m^2+p^3\over (\sqrt {m^2+p^2})^3}>0.
$$
So $p_0$ is unique.
\qed
\begin{lem}\label{selection}  Suppose that $s(x)$ and $k(x)$ are $C^{\infty}$ periodic functions. Given $m\ne 0, n\ne 0, c\in \Rset$ and $a>|m|+\max_{\Bbb T^1}k$, there exists at most one periodic viscosity solution $u$ satisfying
$$
\sqrt {m^2+(n+u')^2}+{c(n+u')s(x)\over \sqrt{m^2+(n+u')^2}}+k(x)=a.
$$
Also, $u\in C^{\infty}(\Bbb T^1)$.
\end{lem}
Proof:  Without loss of generality, we assume that $n>0$. Clearly,  $u\in W^{1,\infty}(\Bbb T^1)$. Due to Lemma \ref{unique}, to prove that it is unique and $C^1$, it suffices to show that
\be\label{positive}
n+u'>0  \quad \text{a.e in $\Bbb T^1$}.
\ee
We claim that $w(x)=nx+u(x)$ is increasing.   In fact, if not, then $w$  must have a local minimum point  since $w(x+1)-w(x)=n>0$.  Suppose $w$ attains local minimum at $x_0$. Then by the definition of viscosity solutions, we have that
$$
|m|+k(x_0)\geq a.
$$
This is a contradiction.  Hence our claim holds. Note that $w'$ can not be zero. Therefore (\ref{positive}) is true.
\qed

Assume that $s,k$ are smooth periodic functions. Fix $m$ and $n$. For $c\geq 0$, we denote $h(c)$ as the unique number such that the following equation
\be\label{generalcell}
\sqrt {m^2+(n+u')^2}+{c(n+u')s(x)\over \sqrt {m^2+(n+u')^2}}+k(x)=h(c)
\ee
admits viscosity solutions. The existence of $h(c)$ follows from \cite{LPV}.
\begin{theo}\label{parameter}  Suppose that $\int_{0}^{1}s(x)=0$ and $s(x)$ is not identically 0.  Also we assume that
 \be\label{match}
 s(x)=0 \quad \text{if $k(x)=\max_{\Bbb T^1}k$.}
 \ee
 If $m$, $n\ne 0$, then\\
(i) $h(c)$ is Lipschitz continuous and $\min_{[0,+\infty)}h=|m|+\max_{\Bbb T^1}k$;\\
(ii) $h'(c)\leq 0$ for a.e $c\in (0,\infty)$;\\
(iii) If $h(c)>|m|+\max_{\Bbb T^1}k$, then $h'(c)<0$;\\
(iv) There exists $c_0$ such that $h(c)=|m|+\max_{\Bbb T^1}k$, for $c\geq c_0$.
\end{theo}
Proof:  (i) and (iv). We claim that
\be\label{Lip}
|h(c_1)-h(c_2)|\leq L|c_1-c_2|
\ee
for $L=\max_{\Bbb T^1}|s|$.  In fact, suppose that $u_1$ and $u_2$ are viscosity solutions of (\ref{generalcell}) for $h(c_1)$ and $h(c_2)$ respectively. We may assume that both of them are $C^1$.  Otherwise, we just apply the routine arguments using the method of ``doubling the number of variables".  Suppose that
$$
u_1(x_0)-u_2(x_0)=\max_{\Bbb T^1}(u_1-u_2).
$$
Then
$$
u_{1}^{'}(x_0)=u_{2}^{'}(x_0)
$$
and
$$
h(c_1)-h(c_2)={(c_1-c_2)s(x_0)(n+u_{1}^{'}(x_0))\over {\sqrt {m^2+(n+u_{1}^{'}(x_0))^2}}}.
$$
Hence (\ref{Lip}) holds.  Now owing to (\ref{match}), it is obvious that
\be\label{lower}
h\geq |m|+\max_{\Bbb T^1}k.
\ee
Since $s$ is not constantly 0,  there exists $x_0\in (0,1)$ such that $\tau=s(x_0)<0$.  Choose $\epsilon>0$ such that $[x_0-\epsilon,x_0+\epsilon]\in (0,1)$ and for $x\in [x_0-\epsilon,x_0+\epsilon]$
$$
s(x)\leq {-\tau\over 2}.
$$
Let $\phi$ be a $C^{\infty}$ periodic function such that $n+\phi'\geq 0$ and within the interval $[0,1]$
$$
\mathrm{supp}(n+\phi')\subset (x_0-\epsilon,x_0+\epsilon).
$$
Now choose $c_0>0$ such that
$$
c_0\tau>2\max_{\Bbb T^1}(n+\phi').
$$
Then it is clear that for all $c\geq c_0$
$$
\max_{\Bbb T^1}\{\sqrt {m^2+(n+\phi')^2}+{c(n+\phi')s(x)\over \sqrt {m^2+(n+\phi')^2}}+k(x)\}=|m|+\max_{\Bbb T^1}k(x).
$$
Suppose that $u(x,c)$ is viscosity solution of equation (\ref{generalcell}) for $c\geq c_0$. Considering the place where $u-\phi$ attains minimum, we derive that for $c\geq c_0$
$$
h(c)\leq |m|+\max_{\Bbb T^1}k(x).
$$
Combining with (\ref{lower}), (i) and (iv) follows.
\medskip

(ii) and (iii): Suppose that $h$ is differentiable at $\bar c$.  If $h(\bar c)=|m|+\max_{\Bbb T^1}k(x)$, then it is clear that $h'(\bar c)=0$.  So let us assume that $h(\bar c)>|m|+\max_{\Bbb T^1}k(x)$.  Since $h$ is continuous, when $c$ is close to $\bar c$,  $h(c)>|m|+\max_{\Bbb T^1}k(x)$. Owing to Lemma \ref{selection},  let $u=u(x,c)$ be the unique smooth solution of equation (\ref{generalcell}) subject to
$\int_{0}^{1}u\,dx=0$. Without loss of generality, we assume that $n>0$. Then $n+u'>0$. Denote
$$
f(p)=\sqrt {m^2+p^2}.
$$
Then equation (\ref{generalcell}) can be rewritten as
$$
f(n+u')+cs(x)f'(n+u')+k(x)=h(c).
$$
Taking derivative with respect to $c$ on both sides, we have that
$$
h'(c)=f'(n+u')u_{c}^{'}+cs(x)f''(n+u')u_{c}^{'}+f'(n+u')s(x).
$$
Note that
$$
f'(n+u')+cs(x)f''(n+u')={m^2(p(x)+cs(x))+p^3(x)\over (\sqrt {m^2+p^2(x)})^{3}}
$$
for $p(x)=n+u'$. Since $p(x)>0$ and $h(c)-|m|>\max_{\Bbb T^1}k$, owing to (\ref{generalcell}), we have that
$$
p(x)+cs(x)>0.
$$
Hence  $f'(n+u')+cs(x)f''(n+u')>0$. Therefore
$$
u_{c}^{'}=-{f's(x)\over f'+f''cs(x)}+{h'(c)\over f'+f''cs(x)}.
$$
Taking integration on both sides over $[0,1]$, we obtain
$$
h'(c)\int_{0}^{1}{1\over f'+f''cs(x)}\,dx=\int_{0}^{1}{s(x)\over 1+a(x)s(x)}\,dx,
$$
where $a(x)={cf''\over f'}>0$.  Since $1+a(x)s(x)>0$, $\int_{0}^{1}s(x)=0$ and $s(x)$ is not identically 0,
$$
\int_{0}^{1}{s(x)\over 1+a(x)s(x)}\,dx<0.
$$
So for $c$ close to $\bar c$
$$
h'(c)<0.
$$
\qed

{\bf \nit Proof of Theorem {\ref{shear}}}. Case 1: If $m=0$,  then it is easy to see that $\lambda=\hat \lambda=|n|$.

Case 2: If  $n=0$,  then $\lambda=\hat \lambda=|m|+\max_{\Bbb T^1}mv$. In fact, let $u$ be periodic viscosity solution of
$$
\sqrt {m^2+(u')^2}+{u'mv'\over \sqrt {m^2+(u')^2}}+mv=\hat \lambda.
$$
Suppose that $u(x_0)=\min_{\Bbb T^1}u$. Then
$$
\hat \lambda\leq |m|+mv(x_0)\leq |m|+\max_{\Bbb T^1}mv.
$$
Combining with (\ref{lowerbound2}), $\hat \lambda=|m|+\max_{\Bbb T^1}mv$.  The proof of the equality $\lambda=|m|+\max_{\Bbb T^1}mv$ is similar.

Case 3: If $mn\ne 0$,  then it follows immediately from Lemma {\ref {parameter}} by choosing $s(x)=mv'$ and $k(x)=mv$.

\section{Comparison of Asymptotic Growth Rate of $s_T$ from G-equation, F-equation and RDA}
\setcounter{equation}{0}
The main purpose of this section is to find quantitative and qualitative similarities between asymptotic growth rate of the turbulent flame speeds predicted by  G-equation, F-equation and RDA.  Theorem {\ref{character}} implies that in two dimensions, these three models predict the bending effect simultaneously. We also provide a simple necessary and sufficient condition of the nonoccurrence of the bending effect. The growth laws of $s_T$  as $a\to +\infty$ predicted by these three models are completely identified for the cat's eye flow and the cellular flow. This partially answers a question posed at the end of \cite{EMS}. Recall that we normalize $d=s_l=\tau_r=1$. Throughout this section, we assume that $V(x)$ is divergence free and has mean zero.  When $n=2$, this is equivalent to saying that there exists a periodic smooth stream function $H$ such that
$$
V(x)=(-H_{x_2},H_{x_1}).
$$
\subsection{A Rough Comparison of $s_T$ from G and F-Equations}
For $A>0$, we scale $V$ to $AV$. For $p\in \Rset^n$, we denote $\alpha_A(p)$  as $s_T$ from the inviscid G-equation which is the effective Hamiltonian of the following cell problem
\be\label{Gcell}
|p+DG|+AV(x)\cdot (p+DG)=\alpha_A(p).
\ee
We also denote $\beta_A(p)$ as the effective Hamiltonian from the following cell problem of invisicd F-equation
\be\label{Fcell}
|p+DF|^2+AV(x)\cdot (p+DF)=\beta_A(p).
\ee
We recall that $\alpha_A(p)$ is the unique number such that equation (\ref{Gcell}) admits approximate periodic viscosity solutions and $\beta_A(p)$ the unique number such that (\ref{Fcell}) admits periodic viscosity solutions.  Note that (\ref{Gcell}) might not have exact solutions due to the lack of coercivity.  Both $\alpha_A(p)$ and $\beta_A(p)$ can be given by inf-max formulas. Precisely speaking,
\be\label{infmaxG}
\alpha_A(p)=\inf_{\phi\in C^1(\Bbb T^n)}\max_{\Bbb T^n}\{|p+D\phi|+AV(x)\cdot (p+D\phi)\}
\ee
and
\be\label{infmaxMS}
\beta_A(p)=\inf_{\phi\in C^1(\Bbb T^n)}\max_{\Bbb T^n}\{|p+D\phi|^2+AV(x)\cdot (p+D\phi)\}
\ee
Hence $\alpha_A(p)$, $\beta_A(p)\leq O(A)$. Since $t^2>t-{1\over 4}$,
\be\label{super}
\beta_A\geq \alpha_A-{1\over 4}.
\ee
We write $\gamma_A(p)$ as $s_T$ from the Majda-Sougandis model along direction $p$ for $\alpha\in (0,1)$. Then by (\ref{MSturbulent}),
$$
\gamma_A(p)=\inf_{\lambda>0}{f'(0)+\beta_A(\lambda p)\over \lambda}.
$$
\begin{lem}\label{middle} Denote $\tau=2\sqrt {f'(0)}$. Then
$$
\tau\alpha_{A\over \tau}(p)\leq \gamma_{A}(p)\leq \beta_A(p)+f'(0).
$$
\end{lem}
Proof:  The right inequality is obvious by choosing $\lambda=1$.  Let us prove the left inequality. Since ${1\over \tau}t^2> t- {\tau\over 4}$, according to the inf-max formulas,
$$
\beta_A(p)\geq \tau\alpha_{A\over \tau}(p)-f'(0).
$$
Combining with the degree 1 homogeneity of $\alpha_A(p)$ with respect to the $p$ variable,  the above lemma holds.
\qed
\begin{rmk}  Note that if we choose $f'(0)\geq {1\over 4}$ as in \cite{EMS},
then $\gamma_A(p)\geq \alpha_A(p)$. This generalizes the computations in \cite{EMS} for
shear flows, which shows that the G-equation model always predicts
slower turbulent flame speeds than the F-equation model. However
if $f'(0)<{1\over 4}$, this
conclusion is not true in general. Here is a simple example based on the shear flow.
Let $n=2$, $V(x_1,x_2)=(v(x_2),0)$ and $p=(1,0)$. Then an easy calculation shows that
$$
\alpha_A(p)=1+A\max_{\Bbb T^1}v
$$
and
$$
\beta_A({\lambda p})=\lambda^2+A\lambda \max_{\Bbb T^1}v.
$$
Then $\gamma_A(p)=A\max_{\Bbb T^1}v+2\sqrt {f'(0)}<\alpha_A(p)$.
\end{rmk}
\subsection{Identical Asymptotic Growth Rate for $G$ and $F$ Equations}
Throughout this subsection, for convenience, we drop the dependence of $\alpha_A(p)$, $\beta_A(p)$ and $\gamma_A(p)$ on $p$.
The following theorem says that $\alpha_A/A$,  $\beta_A/A$ and $\gamma_A/A$ have the same asymptotic limit.
\begin{theo}\label{bending} Given $p\in \Rset^n$.  Denote
$$
c_p=\inf_{\phi\in C^1(\Bbb T^n)}\max_{\Bbb T^n}\{V(x)\cdot (p+D\phi)\}.
$$
Then
$$
\lim_{A\to +\infty}{\alpha_A\over A}=\lim_{A\to +\infty}{\beta_A\over A}=\lim_{A\to +\infty}{\gamma_A\over A}=c_p.
$$
In particular, G-equation and F-equation models predict the bending effect simultaneously.
\end{theo}

\nit Proof: The proof is simple. Owing to the inf-max formula (\ref{infmaxG}),
$$
{\alpha_A\over A}\geq c_p.
$$
Now fix $\epsilon>0$ and choose $\phi_{\epsilon}\in C^1(\Bbb T^n)$ such that
$$
\max_{\Bbb T^n}\{V(x)\cdot (p+D\phi_{\epsilon})\}\leq \epsilon+c_p.
$$
Then
$$
{\alpha_A\over A}\leq {1\over A}\max_{\Bbb T^n}|D\phi_{\epsilon}|+c_p+\epsilon.
$$
Hence
$$
\lim_{A\to +\infty}{\alpha_A\over A}=c_p.
$$
The proof for $\beta_A$ is similar. The proof for $\gamma_A$ follows from Lemma \ref{middle}.
\qed

\nit Now an interesting question is when $c_p>0$ occurs, i.e, there is no bending effect.
To this end, we employ the weak KAM theory to introduce an equivalent formula of $c_p$.

\begin{lem}\label{invariant}
$$
c_p=\max_{\sigma\in \Lambda}\int_{\Bbb T^n}p\cdot V(x)\,d\sigma,
$$
where $\Lambda$ is the collection of all Borel probability measure on $\Bbb T^n$ invariant under the flow $\dot \xi=V(\xi)$.
\end{lem}

\nit Proof: Assume that $\phi\in C^1(\Bbb T^n)$ and $\sigma\in \Lambda$.  We have that
$$
\int_{\Bbb T^n}p\cdot V(x)\,d\sigma=\int_{\Bbb T^n}p\cdot V(x)+D\phi\cdot V\,d\sigma\leq \max_{\Bbb T^2}\{V(x)\cdot (p+D\phi)\}.
$$
Hence
\be\label{oneside}
c_p\geq \max_{\sigma\in \Lambda}\int_{\Bbb T^n}p\cdot V(x)\,d\sigma.
\ee
Denote the Hamiltonian
$$
H_{A}(p,x)={1\over A}|p|^2+V(x)\cdot p.
$$
Then the corresponding Lagrangian is
$$
L_{A}(q,x)={A\over 4}\left |q-V(x) \right |^2.
$$

Write $\bar H_{A}(p)={\beta_A(p)\over A}$. Then $\bar H_{A}(p)$ is a convex function of $p$. Now fix $p$ and  choose $Q_A\in \partial\bar H_{A}(p)$. Here $\partial \bar H_A(p)$ is the set of subdifferentials of $\bar H_A$ at $p$, i.e,
$$
\partial \bar H_A(p)=\{Q\in \Rset^n|\ \bar H(p')\geq \bar H(p)+Q\cdot (p'-p)\  for\ all\ p'\in \Rset^n\}.
$$
Let $\mu_A$ be a Mather measure associated with the Lagrangian $L_A$ with the rotation vector $Q$. Then $\mu_A$ is a Borel  probability Lagrangian flow invariant measure on the phase space $\Rset ^n\times \Bbb T^n$  which minimizes the functional
$$
\int_{\Rset^n\times\Bbb T^n}L(q,x)\,d\mu
$$
among all Borel  probability Lagrangian flow invariant measures $\mu$ on the phase space $\Rset ^n\times \Bbb T^n$ subject to
$$
\int_{\Rset^n\times \Bbb T^n}q\,d\mu_A=Q_A.
$$
See \cite{M} for the existence of Mather measures. Also, suppose that $F$ is a viscosity solution of
$$
{1\over A}|P+DF|^2+V(x)\cdot (p+DF)=\bar H_{A}(p).
$$
Then $F$ is differentiable on the projection of $\mathrm{supp}(\mu)$ to $\Bbb T^n$ and
$$
p+DF(x)={A\over 2}(q-V(x))  \quad \text{for $(q,x)\in \mathrm{supp}(\mu_A)$}.
$$
See \cite{EG}. Due to Lemma \ref{decay},
\be\label{compact}
|q-V(x)|\leq o(1) \quad \text{for $(q,x)\in\mathrm{supp}(\mu_A)$},
\ee
where $\lim_{A\to +\infty}o(1)=0$. Since $\bar H_{A}(0)=0$,  we have that
$$
p\cdot Q_A\geq \bar H_{A}(p)\geq c_p.
$$
Upon a subsequence if necessary, we may assume that
$$
\mu_{A}\rightharpoonup \mu  \quad \text{weakly in $\Rset ^n\times \Bbb T^n$}.
$$
Suppose that $\sigma$ is the projection of $\mu$ on $\Bbb T^n$. Owing to (\ref{compact}),
$$
\begin{array}{ll}
\int_{\Rset ^n\times \Bbb T^n}p\cdot q\,d\mu&=\lim_{A\to +\infty}\int_{\Rset ^n\times \Bbb T^n}p\cdot q\,d\mu_A\\[5mm]
&=\lim_{A\to +\infty}p\cdot Q_A\geq c_p.
\end{array}
$$
and
\be\label{concentration}
\mathrm{supp}\mu\subset \{(q,x)|\  q=V(x)\}.
\ee
Now we show that $\sigma$ is flow invariant. In fact, since $\mu_A$ is Euler-Lagrangian flow invariant, for any $\phi\in C^1(\Bbb T^n)$,
$$
\int_{\Rset ^n\times \Bbb T^n}q\cdot D\phi\,d\mu_A=0.
$$
Sending $A\to +\infty$ and using (\ref{concentration}), we derive that
$$
\int_{\Bbb T^n}D\phi\cdot V\,d\sigma=0.
$$
Therefore $\sigma\in \Lambda$ and
$$
\int_{\Bbb T^n}p\cdot V(x)\,d\sigma\geq c_p.
$$
Combining with  (\ref{oneside}),  Lemma \ref{invariant} holds. \qed

\begin{lem}\label{decay}  Suppose that $F$ is a viscosity solution of (\ref{Fcell}), then
$$
\lim_{A\to +\infty}\sup_{x\in \mathcal{M}}\left |{DF(x)\over A}\right |=0,
$$
where $\mathcal{M}$ is the set where $F$ is differentiable.
\end{lem}

\nit Proof: Throughout this proof, $C$ denotes a constant depending only on $V$ and $|p|$.  Choose $F$ such that $\int_{\Bbb T^n}F\,dx=0$.  Denote $w(x)={p\cdot x+F\over A}$. Then $w$ satisfies that
$$
|Dw|^2+V(x)\cdot w(x)={\beta_A\over A^2}.
$$
Step 1: we claim that
$$
\lim_{A\to +\infty}w=0  \quad \text{uniformly in $\Bbb T^n$}.
$$
In fact,  since $\beta_A\leq O(A)$,  we have that $|Dw|\leq C$. Upon a subsequence if necessary,  we assume that
$$
\lim_{A\to +\infty}w=h \quad \text{uniformly in $\Bbb T^n$}.
$$
Then $h$ is a mean zero periodic viscosity solution of
$$
|Dh|^2+V(x)\cdot Dh=0.
$$
Taking integration on both side,  we obtain that
$$
\int_{\Bbb T^n}|Dh|^2\,dx=0.
$$
Hence $h\equiv 0$.  Our claim holds.

Step 2:  we claim that
$$
\sup_{x\in \mathcal{M}}|Dw(x)|\leq C\left ({1\over \sqrt {A}}+\sqrt{ {\max_{\Bbb T^n}}\left |{F\over A} \right |}\right ).
$$
Choose $x_0\in \Bbb T^n$ such that $w$ is differentiable at $x_0$.  Let $\xi(t): (-\infty,0)$ be the backward characteristics with $\xi (0)=x_0$. Then
\be\label{characteristic}
w(x_0)-w(\xi (-t))=t{\beta_A\over A^2}+{1\over 4}\int_{-t}^{0}|\dot \xi-V(\xi)|^2\,ds.
\ee
Also, $\xi$ satisfies the Euler-Lagrange equation
\be\label{EL}
\ddot \xi=DV(\xi)\cdot \dot \xi-(\dot\xi-V)\cdot DV(\xi)
\ee
and the equality
\be\label{gradient}
Dw(\xi)={\dot \xi-V(\xi)\over 2}.
\ee
Accordingly,  $|\dot \xi|$, $|\ddot \xi|\leq C$.  Recall that $w={p\cdot x+F\over A}$.  Choosing $t=1$ in (\ref{characteristic}), we deduce that
$$
\int_{-1}^{0}|\dot \xi-V(\xi)|^2\,ds\leq {C\over A}(1+\max_{\Bbb T^n}|F|).
$$
So there is a $t_0\in [0,1]$ such that
$$
|\dot \xi(-t_0)-V(\xi(-t_0))|^2\leq  {C\over A}(1+\max_{\Bbb T^n}|F|).
$$
By equation (\ref{EL}),
$$
{d\over ds}|\dot \xi(s)-V(\xi(s))|^2\leq C|\dot \xi(s)-V(\xi(s))|^2.
$$
This implies that
$$
|\dot \xi(0)-V(\xi(0))|^2\leq e^{Ct_0}|\dot \xi(-t_0)-V(\xi(-t_0))|^2.
$$
Combining with (\ref{gradient}),  our claim holds.  Lemma \ref{decay} follows from Step 1 and Step 2.
\qed

\begin{defin} A smooth curve  $\xi: \Rset\to \Rset ^n$ is called an {\it orbit}
if $\dot \xi(t)=V(\xi (t))$.  Moreover, $\xi$ is called {\it periodic orbit} and $T>0$ is called a {\it period} if $\xi$ is an orbit satisfying that
$\xi (T)-\xi (0)$ is an integer vector.  Furthermore,
$$
Q={\xi (T)-\xi (0)\over T}
$$
is called the rotation vector of $\xi$.
\end{defin}

\begin{theo}\label{flowpicture} Given $p\in \Rset ^n$.

(i) For any orbit $\xi$,
\be\label{less}
\limsup_{T\to +\infty}{p\cdot \xi (T)\over T}\leq c_p.
\ee
In particular, if $\xi$ is a periodic orbit with period $T$, then
$$
c_p\geq {p\cdot (\xi (T)-\xi (0))\over T}.
$$
(ii) There exists an orbit $\xi$ such that
$$
\lim_{T\to +\infty}{P\cdot \xi (T)\over T}=c_p.
$$
(iii) When $n=2$, there exists a periodic orbit $\xi$ with period $T$ such that
$$
{p\cdot (\xi (T)-\xi(0))\over T}=c_p.
$$
where $T$ is the period.
\end{theo}

\nit Proof:  (i) In fact, for any $\phi\in C^1(\Bbb T^n)$,
$$
\lim_{T\to +\infty}{1\over T}\int_{0}^{T}D\phi (\xi (t))\cdot V(\xi (t))\,dt=\lim_{T\to +\infty}{\phi(\xi(T))-\phi(\xi (0))\over T}=0.
$$
Hence
$$
\begin{array}{ll}
\limsup_{T\to +\infty}{p\cdot \xi (T)\over T}&=\limsup_{T\to +\infty}{1\over T}{\int_{0}^{t}P\cdot V(\xi (t))\,dt}\\[5mm]
&=\limsup_{T\to +\infty}{1\over T}{\int_{0}^{t}(p+D\phi (\xi(t)))\cdot V(\xi (t))\,dt}\\[5mm]
&\leq \max_{\Bbb T^n}\left\{(p+D\phi)\cdot V\right\}.
\end{array}
$$
So (\ref{less}) holds.

(ii) Choose $\sigma_0\in \Lambda$ such that
$$
\int_{\Bbb T^n}p\cdot V(x)\,d\sigma_0=\max_{\sigma\in \Lambda}\int_{\Bbb T^n}p\cdot V(x)\,d\sigma=c_p.
$$
By the Birkhoff Ergodic Theorem for measure preserving flows, there exists a $\sigma_0$ measurable function $\bar \psi$ such that for $\sigma_0$ a.e $x$
$$
\bar \psi (x)=\lim_{T\to +\infty}{p\cdot \xi_x (T)\over T}
$$
and
$$
\int_{\Bbb T^n}p\cdot V(x)\,d\sigma_0=\int_{\Bbb T^n}\bar \psi(x)\,d\sigma_0,
$$
where $\dot \xi_{x}=V(\xi_{x})$ and $\xi_x(0)=x$. By (\ref{less}),
$$
\bar \psi (x)\leq \int_{\Bbb T^2}p\cdot V(x)\,d\sigma_0.
$$
Hence for $\sigma_0$ a.e $x$
$$
\bar \psi (x)=\int_{\Bbb T^n}p\cdot V(x)\,d\sigma_0.
$$
(iii) Now let us assume $n=2$, according to Poincar\'e recurrence theorem, a.e. in $\sigma_0$, $x$ are recurrent points. Note that if $x$ is recurrent, $\xi_x$ must be periodic. The reason is simple.   Assume that $V(x_1,x_2)=(-H_{x_2}, H_{x_1})$ for some
smooth stream function $H$. Then $H$ is constant along $\xi_{x}$. Hence $\xi_{x}$ must be periodic.

The following is an immediate corollary.

\begin{cor}\label{ns} Assume $n=2$ and $p\in \Rset ^2$. Then
$$
c_p>0
$$
if and only there exists a periodic orbit $\xi$ with period $T$ such that
$$
p\cdot (\xi (T)-\xi (0))>0.
$$
\end{cor}

\begin{lem}\label{vashingdirection} Suppose that $\xi:[0,T]\to \Rset^2$ is a  periodic orbit and $T$ is the period. Denote $Q={\xi (T)-\xi (0)\over T}$. Assume that $Q\ne 0$.  Then\\
(i) $c_p>0$ if $p\cdot Q>0$;\\
(ii) $\alpha_A(p)$, $\beta_A(p)$, $\gamma_A(p)\leq O(1)$ for all $A\geq 0$ if $p\cdot Q=0$.
\end{lem}

\nit Proof:   (i) follows immediately from Corollary \ref{ns}. Let us prove (ii). Throughout this proof, for $x\in \Rset^2$, $\xi_x$ denotes the orbit satisfying $\xi_{x}(0)=x$. Owing to (\ref{super}) and Lemma \ref{middle}, it suffices to prove for $\beta_A$. The strategy is to construct a suitable function  $S\in W^{1,\infty}(\Rset ^2)$ satisfying
$$
V\cdot (p+DS)=0.
$$
Denote $\bar x=\xi (0)$.  Choose a stream function $H$ such that $V=\nabla^\perp H=(-H_{x_2},H_{x_1})$. Since $H$ is constant along $\xi_{x}$ and $DH$ does not vanish, there exists a number
$r>0$ such that for all $x\in B_r(\bar x)$,  $\xi_{x}$ is periodic and $DH$ does not vanish along $\xi_{x}$. Let $\eta: \Rset\to \Rset^2$ satisfy
$$
\dot \eta (t)=DH(\eta (t))
$$
and
$$
\eta(0)=\bar x.
$$
Then there exists $\delta>0$ such that for all $t\in (-\delta, \delta)$\\
(a) $\xi_{\eta (t)}$ is periodic and has nonzero rotation vector parallel with Q \\
(b) $DH(\eta(t))\ne 0$.

Since $H$ is constant along each orbit, (b) implies that
$$
\xi_{\eta (t_1)}\cap \xi_{\eta (t_2)}=\emptyset \quad \text{if $-\delta<t_1<t_2<\delta$}.
$$
and even more
$$
\hat \xi_{\eta (t_1)}\cap \hat \xi_{\eta (t_2)}=\emptyset  \quad \text{if $-\delta<t_1<t_2<\delta$}.
$$
See the end of the introduction section for the definition of $\hat \xi$.

Suppose that
$$
\xi (T)=\bar x+(m,n)
$$
for $(m,n)\in \Zset^2\backslash\{0\}$. Now consider the region $D_0$ which is
 bounded by $\xi$ and $\xi+(-n,m)$. Without loss of generality, we assume that
the positions of $\xi$, $\eta$ and $(m,n)$ are as shown as in Figure 3.
 Let us define $w$ in $D_0$ as follows.
$$
\begin{cases}
w|_{\xi_{\eta(t)}+(-n,m)}=r,\;\;    \quad \text{for $0\leq t\leq {\delta\over 2}$}\\
w|_{\xi_{\eta(t)}+(-n,m)}=r-h(t),\;\;   \quad \text{for ${{\delta}\over 2}\leq t\leq {\delta}$}\\
w=0,\;\;    \quad \text{elsewhere in}\; D_0.
\end{cases}
$$
Here
$$
r=m^2+n^2
$$ and $h$ is a smooth function which is 0 near ${\delta\over 2}$ and $r$ near ${\delta}$.  Clearly, $w$ is a smooth function satisfying
$$
Dw\cdot V\equiv 0   \quad \text{in $D_0$}.
$$
Since $\xi(s)=\xi(s-T)+(m,n)$,  owing the periodicity of $\xi_{\eta(t)}$, we have that as a set
$$
\xi_{\eta(t)}=\xi_{\eta (t)}+(m,n)  \quad \text{for $|t|\leq {\delta}$}
$$
and
$$
D_0=D_0+(m,n).
$$
Hence for $x\in D_0$:
$$
w(x+(m,n))=x.
$$
Now for $i\in \Zset$, we denote
$$
D_i=D_0+i(-n,m).
$$
Then it is clear that $D_i$ is the region bounded by $\xi+i(-n,m)$ and $\xi+(i+1)(-n,m)$.  Also,  $D_i$ and $D_j$ have disjoint interior if $i\ne j$ and
$$
\cup_{i\in \Zset}D_i=\Rset ^2.
$$
We define  $U:\Rset ^2\to \Rset$ as
$$
U(x)=w(x-i(-n,m))+ri,\;\;    \quad \text{if $x\in D_i$}.
$$
Recall that $r=m^2+n^2$. Then $U\in C^{\infty}(\Rset ^2)$ and
$$
DU\cdot V\equiv  0  \quad\text{in $\Rset^2$}.
$$
Let $S(x)=U(x)-x\cdot (-n,m)$. Then it is easy to see that
$$
S(x+(m,n))=S(x+(-n,m))=S(x).
$$
Hence $S$ is $r\Zset^2$ periodic, i.e, for all $v\in \Zset^2$
$$
S(x+rv)=S(x).
$$
Now suppose $p\cdot Q=0$. Note that $Q={(m,n)\over T}$. There exists $\nu\in \Rset$ such that $p=\nu (-n,m)$.  Let $F$ be a viscosity solution of
$$
|p+DF|^2+AV(x)\cdot (p+DF)=\beta_A.
$$
for $p=\nu (-n,m)$. Assume that $F(x_0)-\nu S(x_0)=\min_{\Rset ^2}(F-\nu S)$.  Then
$$
|p+\nu DS(x_0)|^2+AV(x_0)\cdot (p+\nu DS(x_0))\geq \beta_A.
$$
Since $DU\cdot V\equiv 0$, we deduce that
$$
\beta_A(p)\leq |p+\nu DS(x_0)|^2=O(1).
$$
\qed
\begin{figure}
\centering
\includegraphics{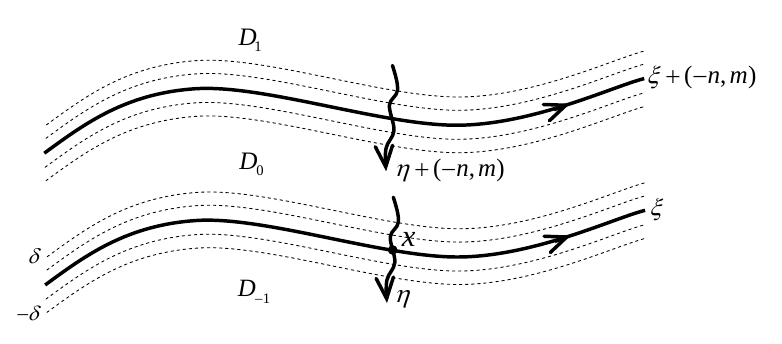}
\caption{Positions of $\xi_{\eta (t)}$.}
\label{fig:flowposition}
\end{figure}
\subsection{Qualitative Similarities  between $c_p$ and the  Asymptotic Growth Rate of $s_T$ from RDA}
Now let us denote $c_{p}^{*}(A)$ as the $s_T(A)$ from the RDA equation (\ref{rd1}),
and given by the variational formula (\ref{rd2}).  According to \cite{Z_10},
$$
\lim_{A\to +\infty}{c_{p}^{*}(A)\over A}=\sup_{w\in \Gamma}{\int_{\Bbb T^n}V\cdot pw^2\,dx\over \int_{\Bbb T^n}w^2\,dx}=c_{p}^{*},
$$
where
$$
\Gamma=\{w\in H^1(\Bbb T^n)|\ V\cdot w=0,\  ||Dw||_{2}^{2}\leq {f'(0)}||w||_{2}^{2}\}.
$$
According to Lemma \ref{invariant}, it is obvious that
$$
c_{p}^{*}\leq c_p.
$$
More interestingly,  Lemma \ref{exchange} says that when $n=2$
\be\label{twolimit}
\lim_{f'(0)\to +\infty}c_{p}^{*}= c_p.
\ee
The equality (\ref{twolimit}) reveals that in 2d, changing the order of two limiting processes
(one sending the flame thickness or reaction time scale to zero,
the other sending the turbulent intensity to infinity) does not change the
asymptotic growth rate of the effective front speed. The following is an analogue of
Lemma \ref{vashingdirection}.
\begin{lem}\label{reactioncircle} Suppose that $\xi:\Rset\to \Rset^2$ is a  periodic orbit with period $T$. Denote the rotation vector $Q={\xi (T)-\xi (0)\over T}$. Assume that $Q\ne 0$.  Then\\
(i) $c_{p}^{*}>0$  if $p\cdot Q>0$;\\
(ii) $c_{p}^{*}(A)\leq O(1)$ for all $A\geq 0$ if $p\cdot Q=0$.
\end{lem}

\nit Proof: (i) Assume that $p\cdot Q>0$. Let $\bar x$, $\eta=\eta(t)$, $\delta$ and $(m,n)$ be the same as in
the proof of Lemma \ref{vashingdirection}, see also the Figure 3.  We may
choose $\delta$ small enough such that for $|t|\leq {\delta}$ and $s\in \Rset$
\be\label{positivedirection}
({\xi_{\eta(t)}(T+s)-\xi_{\eta (t)}(s)})\cdot p>0.
\ee
For $M>0$, we define $w\in W^{1,\infty}(\Bbb T^2)$ as follows
$$
\begin{cases}
w|_{\hat \xi_{\eta(t)}}=M+1-{4t^2\over \delta^2}   \quad \text{for  $|t|\leq {\delta\over 2}$}\\
w=M  \quad \text{elsewhere in $\Bbb T^2$}.
\end{cases}
$$
The definition of $\hat \xi_{\eta (t)}$ is at the end of the introduction section.
Clearly $Dw\cdot V\equiv 0$.  Also, if we choose $M$ large enough,
$$
||Dw||_{2}^{2}\leq {f'(0)}||w||_{2}^{2}.
$$
Denote $C=\cup_{|t|\leq {\delta\over 2}}\hat \xi_{\eta(t)}\subset \Bbb T^2$,
which is a closed invariant set of $\Bbb T^2$.  For all $s>0$,
$$
\int_{C}V(x)\cdot pw^2(x)\,dx=\int_{C}V(\xi_x(s))\cdot pw^2(\xi_x(s))\,dx.
$$
Hence
$$
\begin{array}{ll}
\int_{C}V(x)\cdot pw^2(x)\,dx&={1\over T}\int_{0}^{T}\int_{D}V(\xi_x(s))\cdot pw^2(\xi_x(s))\,dxds\\[5mm]
&={1\over T}\int_{C}\int_{0}^{T}V(\xi_x(s))\cdot pw^2(\xi_x(s))\,dsdx\\[5mm]
&={1\over T}\int_{C}w^2(x)\int_{0}^{T}\dot \xi_x(s)\cdot p\,dsdx\\[5mm]
&=\int_{C}\, w^2(x)({\xi_x(T)-x})\cdot p\,dx
\end{array}
$$
Owing to (\ref{positivedirection}),
$$
\int_{C}V(x)\cdot pw^2(x)\,dx>M^2\int_{C}V\cdot p\,dx.
$$
So according to the definition of $w$
$$
\int_{\Bbb T^2}V(x)\cdot pw^2\,dx>M^2\int_{\Bbb T^2}V(x)\cdot p\,dx=0.
$$
Hence (i) holds.

Now let us prove (ii).  Assume $p=\nu(-n,m)$ for some $\nu\in \Rset$. Due to the formula (\ref{rd2}),
\be\label{VF1}
c_{p}^{*}(A)\leq \kappa_A(p)+f'(0)
\ee
where $\kappa_A(p)$ is the unique number such that the following cell problem has solutions
\be\label{Fe}
-\Delta u+|p+Du|^2+AV(x)\cdot (p+Du)=\kappa_A(p).
\ee
Let $S$ be the smooth $r\Zset^2$ periodic solution from the proof of Theorem \ref{vashingdirection} with $r=m^2+n^2$.  Assume that
$$
u(x_0)-\nu S(x_0)=\min_{\Rset ^2}(u-\nu S).
$$
Combining with the equality that $V(x)\cdot (p+\nu DS)=0$,  we have that
\be\label{VF2}
\kappa_A(p)\leq -\nu\Delta S(x_0)+|p+\nu DS(x_0)|^2=O(1).
\ee
\qed

\begin{lem}\label{rotationvector} Let $\xi:\Rset\to \Rset ^2$ be a periodic orbit with rotation vector $Q$.  Then there must exist another periodic orbit $\tilde \xi$ with rotation vector $\lambda Q$ for some $\lambda<0$.
\end{lem}

\nit Proof:  This is obvious if $Q=0$. So let us assume that $Q\ne 0$.  Since $n=2$, any two periodic orbits must have parallel rotation vectors. So according to the Poincar\'e recurrence theorem, for a.e $x$, $\xi_x$ is a periodic orbit with rotation vector parallel to $Q$.  Hence for a.e $x$,
$$
\lim_{T\to +\infty}{\xi_{x}(T)-x\over T}=\lambda_xQ.
$$
for some $\lambda_x\in \Rset$.  Let $\xi(0)=\bar x$. Then $\lambda_{\bar x}=1$ and when $x$ is close to $\bar x$, $\xi_x$ is also periodic and $\lambda_x>0$. Note that
$$
\begin{array}{ll}
0=\int_{\Bbb T^2}V(x)\cdot Q\,dx&={1\over T}\int_{0}^{T}\int_{\Bbb T^2}V(\xi_x(t))\cdot Q\,dxdt\\[5mm]
&={1\over T}\int_{\Bbb T^2}\int_{0}^{T}V(\xi_x(t))\cdot Q\,dtdx\\[5mm]
&=\int_{\Bbb T^2}{(\xi_x(T)-x)\cdot Q\over T}\,dx\\[5mm]
\end{array}
$$
Sending $T\to +\infty$, by dominant convergence theorem, we derive that
$$
\int_{\Bbb T^2}\lambda_{x}\,dx=0.
$$
So there must exist $\tilde x$ such that $\xi_{\tilde x}$ is periodic and $\lambda_{\tilde x}<0$.
\qed

The following theorem says that  the G-equation model,  the F-equation model and the reaction-diffusion equation model predict the bending effect simultaneously, i.e, if one of these three models predicts the bending effect, so will the other two. It follows immediately from Lemma \ref{vashingdirection}, Lemma \ref{reactioncircle}, (\ref{VF2}) and Lemma \ref{rotationvector}.
\begin{theo}\label{character} Assume $n=2$.  Then one of the following holds.
\medskip

\nit (i)
$$
c_p=c_{p}^{*}=0,  \quad \text{for all $p\in \Rset^2$}.
$$

\nit (ii) There exists a nonzero vector $Q\in \Rset^2$ such that
$$
c_p\geq c_{p}^{*}>0,   \quad \text{for all $p$ such that $p\cdot Q\ne 0$}
$$
and for all $A>0$
$$
\alpha_A(p), \beta_A(p), \gamma_A(p), \kappa_A(p), c_{p}^{*}(A)\leq O(1),   \quad \text{if $|p|$=1 and $p\cdot Q=0$}.
$$
Moreover, case (ii) happens if and only if there exists a periodic orbit with rotation vector $Q\ne 0$.
\end{theo}
\medskip
\nit {For reader's convenience, we recall the notations.}
$$
\begin{cases}
\alpha_A(p)   \quad \text{the turbulent flame speed predicted by the G-equation}\\
\gamma_A(p)   \quad \text{the turbulent flame speed predicted by the F-equation}\\
c_{p}^{*}(A)  \quad \text{the turbulent flame speed prediced by the scalar RDA equation}
\end{cases}
$$
and
$$
\begin{cases}
\beta_A(p)  \quad \text{effective Hamiltonian from the F-equation (\ref{Fcell})}\\
\kappa_A(p)  \quad\text{effective Hamiltonian from the viscous F-equation (\ref{Fe}).}
\end{cases}
$$
Also
$$
c_p=\lim_{A\to +\infty}{\alpha_A(p)\over A}=\lim_{A\to +\infty}{\beta_A(p)\over A}=\lim_{A\to +\infty}{\gamma_A(p)\over A}  \quad \text{(Theorem \ref{bending})}
$$
and
$$
c_{p}^{*}=\lim_{A\to +\infty}{c_{p}^{*}(A)\over A}.
$$

\vspace{5mm}

A good example to demonstrate the above theorem is the following cat's eye flow.

\begin{example}[cat's eye flow]For the cat's eye flow, the stream function is
$H=\sin2\pi x_1 \, \sin2\pi x_2+\, \delta \cos2\pi x_1 \cos2\pi x_2$ for $\delta\in [0,1]$.
Note that when $\delta=0$, it is the cellular flow which belongs to case (i) in Theorem \ref{character}.  For $\delta>0$,
the zero level curve $\{H=0\}$ is a periodic orbit and has a rotation vector parallel to (1,1). So it is the case (ii) in Theorem \ref{character}. Hence
$$
c_p, \; c_{p}^{*}>0,  \quad \text{if $p$ is not parallel to $(-1,1)$}
$$
and for all $A\geq 0$
$$
\alpha_A(p),\; \beta_A(p),\;\gamma_A(p),\; \kappa_A(p),\; c_{p}^{*}(A)\leq O(1), \quad \text{if $p$ is parallel to $(-1,1)$}.
$$
Hence for $\delta\ne 0$, the bending effect only occurs along the direction  parallel to $(-1,1)$.
\end{example}
\subsection{Growth Laws of $s_T$ from Different Models in Cellular Flows}
In the case (i) of Theorem \ref{character},
it remains an important and challenging problem to determine
the exact growth law of $\alpha_A$, $\gamma_A$ and $c_{p}^{*}(A)$ as $A\to +\infty$.
When $V$ is the cellular flow, it is known that
\be\label{loglaw}
\alpha_A=O\left ({A\over \log{A}}\right )   \quad \text{if $p=(1,0)$}
\ee
and
\be\label{quarterlaw1}
c_{p}^{*}(A)=O(A^{1\over 4})  \quad \text{if $p=(1,0)$}.
\ee

Due to (\ref{infmaxG}), (\ref{infmaxMS}), (\ref{infmaxkpp}) and (\ref{rd2}),
it is easy to see that $\alpha_A(p)$, $\beta_A(p)$, $\gamma_A$, $\kappa_A(p)$ and $c_{p}^{*}(A)$ are
all convex as functions of $p$. The symmetry of the stream function $H$ also
implies that these four functions are even and symmetric: $f(p)=f(-p)$ and
$$
f((p_1,p_2))=f((p_2,p_1))=f((-p_1,p_2)).
$$
Note that for $\kappa_A$ and $c_{p}^{*}(p)$, we also need to use the fact that two adjoint compact operators have the same eigenvalues.  Hence (\ref{loglaw}) and (\ref{quarterlaw1}) are true for all $p\ne 0$.  In the following, we show
that $\beta_A$ and $\gamma_A$ also enjoy the ${A\over \log(A)}$ law.  The exact growth law of $\kappa_A$ for cellular flow remains an open problem. See Remark \ref{open} for more details.

\begin{theo}\label{log} If $V$ is the cellular flow, then
$$
\gamma_A(p), \beta_A(p)=O\left ({A\over \log(A)}\right ).
$$
for all $|p|=1$.
\end{theo}

\nit Proof: Throughout this proof, $C$ represents a constant which only depends on $V$. According to Lemma \ref{middle} and the above discussions, it suffices to establish this for $\beta_A(p)$ and $p=(1,0)$. It is not hard to show that for any two points $\bar x$, $\hat x\in Q_1=[0,1]\times [0,1]$, there exists a Lipschitz continuous curve $\xi:[0,L]\to \mathbb {R}^2$ such that $\xi (0)=\bar x$, $\xi (L)=\hat x$ and
$$
|\dot \xi(t)-{\sqrt A}V(\xi (t))|\leq 1  \quad \text{for a.e $t\in [0,L]$}.
$$
Moreover $L\leq C$. Let $U(x)={p\cdot x+F(x)\over \sqrt {A}}$ for $A\geq 1$. Then $U$ satisfies that
$$
|DU|^2+\sqrt {A}V\cdot DU={\beta_A\over A}\leq C.
$$
Hence
$$
|U(x_1)-U(x_0)|\leq C.
$$
So if we choose $F$ such that $\int_{\Bbb T^2}F\,dx=0$,  then
$$
|F|\leq C(\sqrt {A}+1).
$$
Hence due to Step 2 in the proof of Lemma \ref{decay}, we obtain that
\be\label{quarterlaw2}
|p+DF|\leq C({A}^{3\over 4}+1).
\ee
Write $\omega_A=\mathrm{esssup}_{\Bbb T^n}|p+DF|$ and denote $\tilde \alpha_A$ as the effective Hamiltonian of the following modified G-equation
$$
\omega_A|p+D\tilde G|+AV(x)\cdot (p+D\tilde G)=\tilde \alpha_A.
$$
Now we claim that
$$
\beta_A\leq \tilde \alpha_A.
$$
In fact,  suppose that $\phi\in C^1(\Bbb T^n)$ and
$$
\phi (x_0)-F(x_0)=\min_{\Rset ^n}(\phi-F).
$$
Then
$$
|p+D\phi(x_0)|^2+AV(x_0)\cdot (p+D\phi (x_0))\geq \beta_A.
$$
Also, it is easy to see that
$$
|p+\phi (x_0)|\leq \omega_A.
$$
Hence
$$
\omega_A |p+D\phi(x_0)|+AV(x_0)\cdot (p+D\phi (x_0))\geq \beta_A.
$$
So
$$
\max_{\Rset^n}\{\omega_A |p+D\phi|+AV\cdot (p+D\phi)\}\geq \beta_A.
$$
Owing to the inf-max formula (\ref{infmaxG}), my claim holds. Hence
$$
\beta_A\leq \tilde \alpha_A=\omega_A \alpha_{A\over \omega_A}\leq C{A\over \log(A)-\log(\omega_A)}\leq O\left ({A\over \log(A)}\right ).
$$
The last inequality is due to (\ref{quarterlaw2}).

\vspace{5mm}

\nit {\bf Summary:} The growth laws of the turbulent flame speeds from different models in the cellular flow are as follows. For $p\ne 0$,
$$
\begin{cases}
\alpha_A(p)=O\left({A\over \log(A)}\right)\\
\gamma_A(p)=O\left({A\over \log(A)}\right)\\
c_{p}^{*}(A)=O({A^{1\over 4}}).
\end{cases}
$$
A remaining interesting question is what is $\lim_{A\to +\infty}(\gamma_A-\alpha_A)$ or  $\lim_{A\to +\infty}{\gamma_A\over \alpha_A}$  for the cellular flow.  We will investigate this in a future paper.
\section{Viscous G-equation and Viscous F-equation}
\setcounter{equation}{0}
Recall that without loss of generality, we set $d=\tau_r=s_l=1$. For $A>0$, we denote $\chi_A$ as the effective Hamiltonian of the cell problem associated with the viscous G-equation
\be\label{VG}
-\Delta G+|p+DG|+AV(x)\cdot (P+DG)=\chi_A,
\ee
and $\kappa_A$ as the effective Hamiltonian of the following cell problem
\be\label{KPP}
-\Delta S+|p+DS|^2+AV(x)\cdot (p+DS)=\kappa_A.
\ee
Here we omit the dependence of $\chi_A$ and $\kappa_A$ on $p$. As it was mentioned in the introduction part, $\chi_A$ can be viewed as a simplified model of the turbulent flame speed predicted by the curvature G-equation (\ref{curvatureG}). Also, $\kappa_A$ is a simple
upper bound for the turbulent flame speed $s_T$ predicted by the RDA model and the
Majda-Souganidis model for $\alpha=1$. Precisely speaking, owing to (\ref{rd2}),
$$
c_{p}^{*}(A)\leq \kappa_A+f'(0).
$$
The main purpose of this section is to identify the asymptotic growth rate of $\kappa_A$ in two dimensions. Similar to the inviscid case, both $\chi_A$ and $\kappa_A$ are given by inf-max formulas:
\be\label{infmaxvg}
\chi_A=\inf_{\phi\in C^2(\Bbb T^n)}\max_{\Bbb T^n}\{-\Delta \phi+|p+D\phi|+AV(x)\cdot (p+D\phi)\},
\ee
and
\be\label{infmaxkpp}
\kappa_A=\inf_{\phi\in C^2(\Bbb T^n)}\max_{\Bbb T^n}\{-\Delta \phi+|p+D\phi|^2+AV(x)\cdot (p+D\phi)\}.
\ee
Clearly,  $\chi_A$, $\kappa_A\leq O(A)$ and $\chi_A\leq \kappa_A+{1\over 4}$. The following theorem
simply says that the asymptotic behavior in the viscous case can not be larger than that in the inviscid case.
\begin{lem}\label{upperbound}
$$
\limsup_{A\to +\infty}{\chi_A\over A}\leq \limsup_{A\to +\infty}{\kappa_A\over A}\leq c_p.
$$
\end{lem}

\nit Proof:  The first inequality is obvious. Let us prove the second one. For $\epsilon>0$, by routine mollification, there exists a $\phi\in C^2(\Bbb T^n)$ such that
$$
\begin{array}{ll}
\max_{\Bbb T^n}\{V(x)\cdot (p+D\phi_{\epsilon})\}&\leq \inf_{\phi\in C^1(\Bbb T^n)}\max_{\Bbb T^n}\{V(x)\cdot (p+D\phi)\}+\epsilon\\[5mm]
&=c_p+\epsilon.
\end{array}
$$
Then by (\ref{infmaxkpp}),
$$
{\kappa_{A}\over A}\leq  {1\over A}\max_{\Bbb T^n}(-\Delta \phi_{\epsilon}+|p+D\phi_{\epsilon}|^2)+c_p+\epsilon.
$$
So the inequality holds by sending $A\to +\infty$. \qed

Due to the weak nonlinearity in the viscous G-equation, $\lim_{A\to +\infty}{\chi_A}/A$
might be smaller than $c_p$.  See \cite{LXY} for the example of shear flows.

The following theorem says the equality holds for $\kappa_A$ in two
space dimensions. Precisely speaking,
\begin{theo}\label{samelimit}
For $n=2$,
$$
\lim_{A\to +\infty}{\kappa_A\over A}=c_p.
$$
\end{theo}
We first prove several Lemmas.

\begin{lem}\label{smoothinvariant}
$$
\liminf_{A\to +\infty}{\kappa_{A}\over A}\geq \sup_{\sigma\in \hat \Lambda}\int_{\Bbb T^n}p\cdot V\,d\sigma,
$$
where $\hat \Lambda$ is the collection of smooth invariant measures, i.e,
$$
\hat \Lambda=\{\sigma=w^2dx|\ w\in H^1(\Bbb T^n),\  \int_{\Bbb T^2}w^2\,dx=1\ \mathrm{and}\  V(x)\cdot Dw=0\}.
$$
\end{lem}

\nit Proof: Choose $w\in \hat \Lambda$.  Multiply $w^2$ on both sides of equation (\ref{KPP}) and integration by parts.  We deduce that
$$
{1\over A}\int_{\Bbb T^n}wDSDw\,dx+{1\over A}\int_{\Bbb T^n}|p+DS|^2w^2\,dx+\int_{\Bbb T^n}p\cdot V(x)w^2\,dx={\kappa_{A}\over A}.
$$
By Cauchy inequality
$$
|\int_{\Bbb T^n}wDSDw\,dx|\leq {1\over 2}( \int_{\Bbb T^n}|DS|^2w^2\,dx+\int_{\Bbb T^n}|Dw|^2\,dx).
$$
Sending $A\to +\infty$,  we have that
$$
\liminf_{A\to +\infty}{\kappa_{A}\over A}\geq \int_{\Bbb T^n}p\cdot V(x)w^2\,dx.
$$
\qed

Recall that Lemma \ref{invariant} says that
$$
c_p=\max_{\sigma\in \Lambda}\int_{\Bbb T^n}p\cdot V(x)\,d\sigma,
$$
where $\Lambda$ is the collection of all Borel probability measure on $\Bbb T^n$ invariant under the flow $\dot \xi=V(\xi)$. Combining Lemma \ref{upperbound} and \ref{smoothinvariant}, the following is an immediate corollary.

\begin{cor} If
\be\label{dense}
\sup_{\sigma\in \hat \Lambda}\int_{\Bbb T^n}p\cdot V(x)\,d\sigma=\max_{\sigma\in \Lambda}\int_{\Bbb T^n}p\cdot V(x)\,d\sigma.
\ee
Then $\lim_{A\to +\infty}{\kappa_A\over A}=c_p$.
\end{cor}

As we mentioned after (\ref{twolimit}) in section 3.2, the physical meaning of (\ref{dense}) is that changing
the order of limits of sending the flame thickness (or reaction time scale) to
zero and sending the turbulent intensity to infinity does not change the front speed asymptotic growth rate.
It is not clear to us when (\ref{dense}) holds.  Nevertheless we are able to establish it in two dimensions.
\begin{lem}\label{exchange} When $n=2$
\be\label{dense1}
\sup_{\sigma\in \hat \Lambda}\int_{\Bbb T^n}p\cdot V(x)\,d\sigma=\max_{\sigma\in \Lambda}\int_{\Bbb T^n}p\cdot V(x)\,d\sigma.
\ee
\end{lem}

\nit Proof: If $c_p=0$,  then (\ref{dense1}) is obvious.  So let us assume that $c_p>0$.
Then due to Corollary {\ref{ns}, there exists a periodic orbit $\xi_{\bar x}: [0,T]\to \Rset ^2$ such that $\xi_{\bar x}(0)=\bar x$ and
$$
c_p={p\cdot (\xi_{\bar x}(T)-\bar x)\over T}={1\over T}\int_{0}^{T}V(\xi (t))\cdot p\,dt=\int_{\Bbb T^2}p\cdot V(x)\,d\sigma
$$
where $T$ is the period and $\sigma$ the Borel probability invariant measure defined as
 $$
 \int_{\Bbb T^2}f(x)\,d\sigma={1\over T}\int_{0}^{T}f(\xi (t))\,dt  \quad \text{for all $f\in C(\Bbb T^2)$}.
 $$
 The strategy of the proof is to construct a sequence of smooth invariant measures to approximate $\sigma$ which implies (\ref{dense1}). Let $\delta$, $\eta=\eta(t)$ be from the proof of Lemma \ref{vashingdirection}. See also Figure 3.  Choose a sequence of functions $\{\phi_n\}$ such that for all $n\in\Nset$,
 $\phi_n\in C^{\infty}(\Bbb T^2)$ and\\
(i)
$$
\phi_{n}|_{\hat \xi}=1;
$$
(ii)
$$
\phi_{n}|_{\hat \xi_{\eta(t)}}=c_t   \quad \text{for some $c_t\in \Rset$};
$$
(iii)
$$
\phi_{n}=0  \quad \text{in $\Bbb T^2\backslash {\cup_{|t|\leq {\delta\over n}}\hat \xi_{\eta(t)}}$};
$$
(iv)
$$
\int_{\Bbb T^2}\phi_{n}^{2}\,dx=1.
$$
Then it is clear that
$$
V(x)\cdot D\phi_n=0.
$$
Let us assume that
$$
\phi_{n}^{2}\rightharpoonup \hat \sigma  \quad \text{weekly as Random measures on $\Bbb T^n$}.
$$
Then $\hat \sigma$ must an invariant probability measure supported on $\hat \xi$. Hence
$$
\hat \sigma=\sigma.
$$
\qed

\medskip
\nit {\bf Proof of Theorem \ref{samelimit}:} It follows immediately from Corollary 4.1 and Lemma \ref{exchange}.
\begin{rmk}\label{open} It remains an open problem whether $\lim_{A\to +\infty}{\kappa_A\over A}$
exists for $n\geq 3$. We also do not know the exact growth law
of $\kappa_A$ in cellular flow except for the sublinearity. Since $c_{p}^{*}(A)\leq \kappa_A+f'(0)$,
$\kappa_A$ has a growth lower bound of $O(A^{1\over 4})$. However,
it is proved in \cite{LXY_11} that $\chi_A$ is uniformly bounded as $A\to +\infty$, i.e,  for all $A\geq 0$
$$
{\chi_A(p)\over |p|}\leq  O(1).
$$
\end{rmk}
A remaining challenging and very interesting question is whether the original curvature equation (\ref{curvatureG}) will predict similar uniform boundness of the turbulent flame speed.  The physical meaning is whether the curvature effect could have essential contributions to the ``strong bending effect" observed in experiments \cite{R}.
\medskip
\section{Appendix: Review of Viscosity Solutions of \\Hamilton-Jacobi Equation and Cell Problem.}
\setcounter{equation}{0}
Suppose that $\Omega$ is an open set in $\Rset^n$ and $H=H(p,x)\in C(\Rset^n\times \Omega)$.  Let us look at the following Hamilton-Jacobi equation
\be\label{HJ}
H(Du,x)=0  \quad\text{in $\Omega$}.
\ee
It is well known that  the above equation might not have classical ($C^1$) solutions.  Solutions need to be interpreted in the viscosity sense. See the {\it User's guide} \cite{CIL} for reference.
\medskip

\nit{\bf$\bullet$ Viscosity subsolution (supersolution):} We say that $u\in C(\Omega)$ is a {\it viscosity subsolution (supersolution)} of
$$
H(Du,x)=0  \quad \text{in $\Omega$}
$$
if for $\phi\in C^1(\Omega)$ and $x_0\in \Omega$ satisfying that
$$
\phi (x_0)-u(x_0)\leq (\geq)\phi (x)-u(x)   \quad \text{for $x\in \Omega$}
$$
then
$$
H(\phi (x_0),x_0)\leq (\geq)0.
$$

\medskip

\nit{\bf$\bullet$ Viscosity solution:} We say that $u\in C(\Omega)$ is a {\it viscosity solution} if it is both a viscosity subsolution and supersolution.

\medskip

\nit{\bf$\bullet$ Approximate viscosity solution:}  We say that the equation
$$
H(Du,x)=0
$$
admits approximate viscosity solutions if for any $\epsilon>0$, there exists $u_{\epsilon}\in C(\Omega)$ such that it is both a viscosity subsolution of
$$
H(Du_{\epsilon},x)=\epsilon
$$
and a viscosity supersolution of
$$
H(Du_{\epsilon},x)=-\epsilon.
$$

\medskip

\nit{\bf $\bullet$ Connection between viscosity solutions and classical solutions:}  Any $C^1$ solution of equation (\ref{HJ}) must be a viscosity solution.   If $u$ is a viscosity solution and is differentiable at $\hat x\in \Omega$, then it satisfies the equation at $\hat x$, i.e,
$$
H(Du(\hat x), \hat x)=0.
$$

\medskip

\nit{\bf$\bullet$  Coercivity of $H$ and Lipschitz continuity of viscosity solutions}:  Suppose that $H$ is coercive in the $p$ variable, i.e
$$
\lim_{|p|\to +\infty}H(p,x)=+\infty  \quad \text{uniformly in $\bar \Omega$},
$$
then any viscosity solution of equation (\ref{HJ}) is $W^{1,\infty}$.

\medskip

\nit{\bf$\bullet$ Cell problem, Homogenization and the effective Hamiltonian}:  Suppose that $H\in C(\Rset^n\times \Bbb T^n)$ and is coercive in the $p$ variable.  It was proved in \cite{LPV} that for any $P\in \Rset^n$, there exists a UNIQUE number $\bar H(P)$ such that the following equation
$$
H(P+Du,x)=\bar H(P)  \quad \text{in $\Rset ^n$}
$$
admits periodic viscosity solutions.  This is the so called ``cell problem" in the homogenization theory.
The quantity $\bar H(P)$ is called {\it effective Hamiltonian}.  For $\epsilon>0$, assume that $u^{\epsilon}$ is the unique viscosity solution of
$$
\begin{cases}
u_{t}^{\epsilon}+H(Du^{\epsilon},{x\over \epsilon})=0\\
u^{\epsilon}(x,0)=g(x).
\end{cases}
$$
Then as $\epsilon\to 0$, $u^{\epsilon}$ uniformly converge to the $\bar u$ which is the unique solution of the following effective equation
$$
\begin{cases}
\bar u_{t}+\bar H(D\bar u)=0\\
\bar u(x,0)=g(x).
\end{cases}
$$
The cell problem and homogenizaton analysis were performed for fully nonlinear PDEs in \cite{Evans}.

\medskip

\nit{\bf$\bullet$  Inf-max formula for convex effective Hamiltonian:}  If $H$ is also convex in the $p$ variable, the $\bar H(P)$ is given by the following inf-max formula. See \cite{Gomes} for instance
$$
\bar H(P)=\inf_{\phi\in C^1(\Bbb T^n)}\max_{\Bbb T^n}H(P+D\phi,x).
$$

\section{Concluding Remarks}
We analyzed the turbulent flame speeds ($s_T$) of
G-equation, the analogous quadratically nonlinear Hamilton-Jacobi equation (the F-equation) and the reaction-diffusion-advection equation in
the presence of steady and periodic compressible flows, shear flows, and
incompressible flows (cellular and cat's eye flows). The strain effect slows down flame propagation in 1d compressible flow and shear flows. If $f'(0)$ is no less than a threshold value, the G-equation always predicts smaller turbulent flame speeds than the F-equation model. The F and G equations share the same asymptotic growth rates in incompressible flows, in particular,
the $s_T$'s obey the same sublinear growth laws ${A\over \log(A)}$ in cellular flows. Moreover,
in two dimensions, the three models predict the bending effect simultaneously. Moreover, in 2d,
changing the order of limits of sending the flame thickness (or reaction time scale) to zero and
sending the turbulent intensity to infinity does not change the asymptotic growth rate.
The nonoccurrence of the  bending effect is characterized precisely by the existence of
periodic orbits on the torus and their rotation vectors in the underlying advective flow field.
The effective Hamiltonian of the viscous F-equation has the same growth rate as the
inviscid F and G-equations in two dimensions. Asymptotic growth rate of
 the effective Hamiltonian from the viscous F-equation in higher dimensions will be explored further in a future work.

\bibliographystyle{plain}

\end{document}